\documentclass[11pt]{article}
\usepackage[dvips]{graphics}
\input{amssym.def}
\input{amssym.tex}

\title{\Large\bf  Inverse antiplane problem on $n$ uniformly stressed inclusions }
\author{\bf {\sc Y.A. Antipov}\\ 
Department of Mathematics, Louisiana State University\\
Baton Rouge LA 70803, USA}

\date{}

\newcommand{\I}{\mathop{\rm Im}\nolimits}
\newcommand{\R}{
\mathop{\rm Re}\nolimits}
\newcommand{\const}{\mbox{const}}

\newcommand{\Md}{\partial}

\newcommand{\ov}[1]{\overline{#1}}

\newcommand{\Ga}{\alpha}
\newcommand{\Gb}{\beta}
\newcommand{\Gd}{\delta}

\newcommand{\Gf}{\phi}
\newcommand{\Gvf}{\varphi}
\newcommand{\Gg}{\gamma}

\newcommand{\Gk}{\kappa}

\newcommand{\Gl}{\lambda}
\newcommand{\Gn}{\eta}

\newcommand{\Gr}{\rho}

\newcommand{\Gs}{\sigma}

\newcommand{\Go}{\omega}
\newcommand{\Gx}{\xi}

\newcommand{\Gz}{\zeta}
\newcommand{\GD}{\Delta}
\newcommand{\GF}{\Phi}

\newcommand{\GL}{\Lambda}

\newcommand{\CC}{{\cal C}}
\newcommand{\CD}{{\cal D}}

\newcommand{\CR}{{\cal R}}

\newcommand{\beq}{\begin{equation}}
\newcommand{\eeq}{\end{equation}}
\newcommand{\barr}{\begin{eqnarray}}
\newcommand{\earr}{\end{eqnarray}}
\newcommand{\beqn}{\begin{equation*}}
\newcommand{\eeqn}{\end{equation*}}
\newcommand{\barrn}{\begin{eqnarray*}}
\newcommand{\earrn}{\end{eqnarray*}}

\newcommand{\fr}{\frac}

\begin{document}
\maketitle

\begin{abstract}

The inverse problem of antiplane elasticity on determination of the profiles of $n$
uniformly stressed inclusions is studied. The inclusions are in ideal contact with the surrounding matrix, the stress field
inside the inclusions is uniform, and at infinity the body is subjected to antiplane uniform shear.
 The exterior of the inclusions, an $n$-connected domain, is treated as the image by a 
 conformal map of an $n$-connected slit domain with the slits lying in the same line.
 The inverse problem is solved by quadratures by reducing it to two Riemann-Hilbert problems
 on a Riemann surface of genus $n-1$. Samples of two and three symmetric and non-symmetric 
uniformly stressed inclusions are reported.

 \end{abstract}

\setcounter{equation}{0}

\section{Introduction}

Inverse boundary value problems for partial differential equations have been of interest  for many researchers since the work by Riabuchinsky (1929) on determination of the boundary of a domain if the function is harmonic inside
and its values and normal derivative are known on the boundary.  An extensive survey of early results on 
inverse boundary value  problems and their applications to the theory of filtration and hydroaerodynamics
was given by  Aksent'ev et al (1980).
Inverse problems of elasticity on determination of the shapes of curvilinear  cavities and inclusions with elastic constants different from those of the surrounding matrix when the fields inside the inclusions have
prescribed properties have been attracting applied mathematicians, mechanical engineers, and
materials scientists 
since the work by Eshelby (1957).
He showed that if the unbounded elastic body 
is uniformly loaded at infinity, and the body has an elliptic or ellipsoidal inclusion with different elastic constants, 
then the  stress field is uniform inside the inclusion. Eshelby conjected that there do not exist other shapes of 
a single inclusion with such a property.  Sendeckyj (1970) proved this conjecture
in the plane and anti-plane cases. An alternative proof for the antiplane
case by the method of conformal mappings was proposed by Ru and Schiavone (1996). 
The model problem on a plane of finite thickness reinforced by a single inclusion of circular
cross-section and subjected to far-field anti-plane shear (mode III) loading (as well as modes I/II loading) was considered 
by Chaudhuri (2003). 

Cherepanov  (1974) studied  the inverse elasticity problem for a plane uniformly
loaded at infinity and having $n$ holes. The holes boundaries are  
subjected to constant normal and tangential traction; these boundaries are not prescribed and have to be 
determined from the condition that the tangential normal stress $\Gs_t=\Gs=\const$ in all the 
contours.  Cherepanov employed the method of conformal mappings to
transform an  $n$-connected slit domain into the elastic domain and reduced the problem to two homogeneous Schwarz problems of the theory of analytic functions on  $n$ slits. These
problems were solved 
(Cherepanov, 1974)  for the symmetric case of two holes.
To solve the Cherepanov problem for any $n$-connected domain, Vigdergauz (1976) proposed
to employ a circular map from the exterior of $n$-circles onto the $n$-connected elastic domain,
integral equations, and the method of least squares for their numerical solution.
An explicit representation in terms of the Weierstrass elliptic function for the profile
of an inclusion in the case of a doubly periodic structure was given by Grabovsky and Kohn (1995).
Recently,  Antipov (2017) developed a method of the Riemann-Hilbert problem on a Riemann
surface of genus $n-1$ to reconstruct  a family of conformal maps for the Cherepanov
problem in the general case of two and three equal-strength cavities and for the case $n\ge 4$ when the preimage
of the plane with holes is a slit domain with the slits lying in the real axis. 

Kang et al (2008) 
analyzed the case of  two antiplane inclusions with the Eshelby uniformity property.
They employed the Weierstrass zeta function and the Schwarz-Christoffel formula to determine
the profile of two symmetric inclusions. 
A method of Laurent series and a conformal mapping from an annulus to a doubly connected domain
 to reconstruct the shape of two inclusions with uniform stresses was
 applied by Wang (2012). Other numerical approaches were applied by Liu (2008) and Dai et al (2017). 

In this paper we aim to propose an exact method of conformal mappings and the Riemann-Hilbert problem
on a Riemann surface of genus $n-1$ for the inverse antiplane problem on  $n$ 
 inclusions. The inclusions may have different shear moduli and are in ideal contact with the surrounding elastic matrix 
 subjected at infinity to uniform antiplane shear $\tau_{13}=\tau_1^\infty$
and  $\tau_{23}=\tau_2^\infty$. The profiles of the inclusions are not prescribed and have to be determined
from the condition that the stress field inside all the inclusions is uniform, $\tau_{13}=\tau_1$
and  $\tau_{23}=\tau_2$. 

In section 2, we assume that such inclusions exist and treat the $n$-connected exterior of the inclusions,
$D^e$, as the image of a slit domain $\CD^e$ with the slits lying in the real axis. 
We reduce the problem of  determination of the conformal maps to two inhomogeneous Schwarz problems
to be solved consecutively. We analyze the particular case $n=1$
  in section 3 and show that the profile of the inclusion is an ellipse by employing two
  maps, a circular and a slit conformal mappings.  
Section 4 treats the case $n=2$, solves the two Schwarz problems by
reducing them to two Riemann-Hilbert problems on an elliptic  Riemann surface and
derives the associated conformal map in terms of singular and elliptic integrals.
In section 5  we analyze the case $n=3$ and derive the conformal mapping  in terms of singular and genus-2
hyperelliptic integrals.
A class of $n$-connected domains $D_*^e$  which may be considered as images by a conformal map
of $n$ slits lying 
in the same line  is considered in section 6 (for $n=1,2,3$, each domain $D^e$ belongs to this class).  A closed-form representation in terms of hyperelliptic integrals
 for a family of such conformal mappings is constructed
 by employing the theory of the Riemann-Hilbert problem on a genus-($n-1$)
Riemann surface.

\vspace{.1in}

\setcounter{equation}{0}

\section{Formulation}\label{s2}

Consider the following inverse problem of  elasticity:

{\sl Suppose an infinite isotropic solid contains $n$ curvilinear inclusions $D_0$, $D_1, \ldots D_{n-1}$. Let
the shear moduli of the inclusions $D_j$ and the matrix $D^e={\Bbb R}^2\setminus D$,
$(D=\cup_{j=0}^{n-1}D_j)$
be $\mu_j$ and $\mu$, respectively. It is assumed that the inclusions are in ideal contact  with 
the matrix, and the solid $D^e\cup D$ is in a state of antiplane shear  due 
to constant shear stresses applied at infinity, $\tau_{13}=\tau_1^\infty$, $\tau_{23}=\tau_2^\infty$. 
It is  required to determine the boundaries of the inclusions, $L_j$, such that the stresses
$\tau_{13}$ and $\tau_{23}$ are constant in all the inclusions $D_j$, $\tau_{13}=\tau_1$, $\tau_{23}=\tau_2$,
$j=0,1,\ldots,n-1$. }

\vspace{2mm}

Let $u$ be the $x_3$-component of the displacement vector. Then $\tau_{j3}=\mu\Md u/\Md x_j$
($j=1,2$), $(x_1,x_2)\in D^e$, and 
\beq
u\sim \mu^{-1}(\tau_1^\infty x_1+\tau_2^\infty x_2)+\const,\quad  x_1^2+x_2^2\to\infty.
\label{2.1}
\eeq
Since the stresses $\tau_{12}$ and $\tau_{13}$ are constant in the inclusions, the  $x_3$-displacements
$u_j$ for ($x_1,x_2)\in D_j$ are linear functions
\beq
u_j=\mu_j^{-1}(\tau_1 x_1+\tau_2 x_2)+d_j', \quad  (x_1,x_2)\in D_j, \quad j=0,1,\ldots,n-1,
\label{2.2}
\eeq
and $d_j'$ are real constants. The boundary conditions of ideal contact imply
that the traction component $\tau_{\nu3}$ and  the $x_3$-component of the displacement 
are continuous  through the contours $L_j$,
\beq
\mu\fr{\Md u}{\Md \nu}=\mu_j\fr{\Md u_j}{\Md \nu}, \quad u=u_j, \quad (x_1,x_2)\in L_j, \quad j=0,1,\ldots,n-1,
\label{2.3}
\eeq
where $\fr{\Md}{\Md \nu}$ is the normal derivative.

The displacement $u$ and $u_j$ are harmonic functions in the domains $D^e$ and $D_j$, respectively. 
Let $v$ and $v_j$ be their harmonic conjugates and denote $z=x_1+ix_2$.
Then the functions
$\Gf(z)=u(x_1,x_2)+iv(x_1,x_2)$ and $\Gf_j(z)=u_j(x_1,x_2)+iv_j(x_1,x_2)$ are analytic in the corresponding domains.
In terms of these functions the ideal contact boundary conditions can be written as
\beq
\fr{\Gk_j+1}{2}\Gf_j(z)-\fr{\Gk_j-1}{2}\ov{\Gf_j(z)}=\Gf(z)+ib_j, \quad z\in L_j, \quad j=0,1,\ldots,n-1,
\label{2.4}
\eeq
where $\Gk_j=\mu_j/\mu$, and $b_j$ are real constants. To verify the equivalence of the boundary conditions (\ref{2.3})
 and (\ref{2.4}), one may use the Cauchy-Riemann conditions $\fr{\Md u}{\Md\nu}=\fr{\Md u}{\Md s}$,
where $\fr{\Md}{\Md s}$ is the tangential derivative. Since the functions $u_j$ are known and given by (\ref{2.2}),
due to the Cauchy-Riemann conditions they are defined up to arbitrary constants
by
\beq
\Gf_j(z)=\fr{\bar\tau z}{\mu_j}+d_j, \quad z\in D_j, \quad j=0,1,\ldots,n-1,
\label{2.5}
\eeq
where $\bar\tau=\tau_1-i\tau_2$, $d_j=d_j'+id_j''$, $d_j''$ are real constants.
In view of the relations  (\ref{2.5}), instead of the function $\Gf(z)$, it is convenient to deal with the function 
$f(z)=\Gf(z)-\bar\tau z/\mu$, $z\in D^e$.
The new function $f(z)$ is analytic in the domain $D^e$, satisfies the boundary condition
\beq
f(z)=\fr{1}{\Gl_j}\R[\bar\tau z]+d_j'+ia_j, \quad z\in L_j, \quad j=0,1,\ldots,n-1,
\label{2.6}
\eeq
and, since $\Gf(z)\sim\bar\tau^\infty z/\mu +\const$, $z\to\infty$,  the condition at infinity
\beq
f(z)\sim \fr{(\bar\tau^\infty-\bar\tau)z}{\mu}+\const, \quad z\to\infty.
\label{2.7}
\eeq
Here, $\Gl_j=\mu_j/(1-\Gk_j)$, $a_j=\Gk_j d_j''-b_j$ are real constants, and $\bar\tau^\infty=\tau^\infty_1-i\tau^\infty_2$.

Let $z=\Go(\Gz)$ be a conformal map that transforms an $n$-connected canonical domain $\CD^e$
into the physical domain $D^e$. An example of the domain $\CD^e$ is the exterior of $n$ circles $\CC_j$ ($j=0,1,\ldots,n-1$).  By scaling and rotation, it is always possible to achieve that one of
the circles say, $\CC_0$, is of unit radius and centered at the origin and, in addition, the center of another
circle say, $\CC_1$, falls on the real axis. If the circular map meets the condition $\Go(\infty)=\infty$, then the radius of $\CC_1$, the complex centers and the radii 
of the rest $n-2$ circles cannot be selected arbitrarily. If this map is chosen, the problem of 
determining the function $f(\Go(\Gz))$ and the map $\Go(\Gz)$ itself 
reduces to two Schwarz problem for the circular domain ${\Bbb C}\setminus \cup_{j=0}^{n-1} K_j$
with $\Md K_j=\CC_j$. It may be solved by reducing the problem to two Riemann-Hilbert problems of the theory
of automorphic functions (Chibrikova \& Silvestrov, 1978; Mityushev \& Rogosin, 2000; Antipov \& Silvestrov, 2007;
Antipov \& Crowdy,  2007).

Another choice of the canonical domain  $\CD^e$ is the exterior of $n$ slits $l_j$, $j=0,1,\ldots,n-1$, 
such that $L_j$ are the images of the slits $l_j$ (Antipov \& Silvestrov, 2007; Antipov, 2017). The advantage of a slit map over a circular map is that
the solution is delivered by quadratures, not in a series form. The disadvantage is that although the slit map method
works for any simply-, doubly-, and triply-connected domains, 
for $n\ge 4$, the method is applicable for a certain subclass of domains only, namely when all the slits 
lie in the same line.

It is known  (Keldysh, 1939; Courant, 1950)  that there exists an analytic function $z=\Go(\Gz)$
that conformally maps the extended complex $\Gz$-plane ${\Bbb C}\cup\infty$ cut along $n$ segments
parallel to the real $\Gz$-axis onto the $n$-connected domain $D^e$ in the $z$-plane. For this map,
the infinite point $z=\infty$ is the image of a certain
point $\Gz=\Gz_\infty\in\CD^e$, and in the vicinity of that point the conformal map $\Go(\Gz)$ can be represented as
\beq
\Go(\Gz)=\fr{c_{-1}}{\Gz-\Gz_\infty}+c_0+\sum_{j=1}^\infty c_j(\Gz-\Gz_\infty)^j
\label{2.8}
\eeq
if $\Gz_\infty$ is a finite point, and 
\beq
\Go(\Gz)=c_{-1}\Gz+c_0+\sum_{j=1}^\infty \fr{c_j}{\Gz}
\label{2.9}
\eeq
otherwise. Here, $c_{-1}=c_{-1}'+ic_{-1}''$.

Denote $f(\Go(\Gz))=F(\Gz)$. From the boundary condition (\ref{2.6}) regardless of which map is employed 
we deduce  that the functions $F(\Gz)$
and $\Go(\Gz)$ satisfy the following two Schwarz problems to be solved consecutively.

{\it Find two functions $F(\Gz)$ and $\Go(\Gz)$ analytic in the domain $\CD^e$ and continuous up to the boundary
$l=\cup_{j=0}^{n-1} l_j$ such that
\beq
\I F(\Gz)=a_j, \quad \Gz\in l_j, \quad j=0,1\ldots,n-1,
\label{2.10}
\eeq
and
\beq
\R [\bar\tau\Go(\Gz)]=\Gl_j[\R F(\Gz)-d_j'], \quad \Gz\in l_j, \quad j=0,1\ldots,n-1.
\label{2.11}
\eeq
In the vicinity of the point $\Gz_\infty$, both of the  functions have a simple pole. 
If the point $\Gz_\infty$ is a finite point and $\Go(\Gz)$ admits the representations
(\ref{2.8}), then 
\beq
F(\Gz)\sim \fr{\bar\tau^\infty-\bar\tau}{\mu}\fr{c_{-1}}{\Gz-\Gz_\infty}, \quad \Gz\to\Gz_\infty.
\label{2.12}
\eeq
Otherwise, if $\Gz_\infty=\infty$ and (\ref{2.9}) holds, then 
\beq
F(\Gz)\sim \fr{\bar\tau^\infty-\bar\tau}{\mu}c_{-1}\Gz, \quad \Gz\to\infty.
\label{2.13}
\eeq
In addition, the function $\Go: l_j\to L_j$ ($j=0,\ldots, n-1$) has to be univalent, and the interiors of the images of the contours $l_j$, the domains $D_j$, are disjoint sets.
}

In what follows we apply both maps in the case $n=1$ and proceed with the second choice of the parametric 
domain $\CD^e$ by employing a slit map to
restore the contours $L_j$ ($j=0,1,\ldots,n-1$) when $n\ge 2$.

\setcounter{equation}{0}
\section{Single inclusion}\label{s3}

We start with the simplest case of a single inclusion and derive the solution by exploiting
 a circular and a slit maps.

\subsection{$n=1$: a circular map}\label{s3.1}

Without loss of generality $\CC_0$ is the unit circle centered at the origin, $\Gz_\infty=\infty$, and $a_0=0$.
The solution of the Schwarz problem (\ref{2.10}), (\ref{2.13}) for the unit circle $\CC_0$ 
is given by
\beq
F(\Gz)=\Gb_0-i\Gb_1\Gz+i\bar\Gb_1\Gz^{-1},
\label{3.1}
\eeq
where
\beq
\Gb_1=\Gb_1'+i\Gb_1''=\fr{\bar\tau^\infty-\bar\tau}{\mu}ic_{-1},
\label{3.2}
\eeq
and  $c_{-1}$ and $\Gb_0$ are real constants. The solution of the second Schwarz problem (\ref{2.11}), (\ref{2.9}) can be represented in the form
\beq
\bar\tau\Go(\Gz)=\Gg_{-1}\Gz^{-1}+\Gg_0+\Gg_1\Gz, 
\label{3.3}
\eeq
where $\Gg_j=\Gg_j'+i\Gg_j''$, $j=-1,0,1$. On substituting the expressions (\ref{3.3}) and (\ref{3.1}) into
the boundary condition (\ref{2.11}) and replacing $\Gz$ by $e^{i\Gvf}$, $0\le\Gvf\le2\pi$, we derive
\beq
\Gg_0'=(\Gb_0-d_0')\fr{\mu_0}{1-\Gk_0}, \quad
\Gg_{-1}'+\Gg_1'=\fr{2\Gb_1''\mu_0}{1-\Gk_0}, \quad
\Gg_{-1}''-\Gg_1''=\fr{2\Gb_1'\mu_0}{1-\Gk_0}.
\label{3.4}
\eeq
Finally, by using formula (\ref{2.9})  and the relations (\ref{3.4})
we determine the function $\Go(\Gz)$ up to an additive complex constant $\Gg$
\beq
\Go(\Gz)=c_{-1}\left(\Gz+\fr{\Gd}{\Gz}\right)+\Gg,
\label{3.5}
\eeq
where
\beq
\Gd=\fr{2\Gk_0\tau^\infty-(\Gk_0+1)\tau}{(1-\Gk_0)\bar\tau}.
\label{3.6}
\eeq
Let $\Gk_0\ne 1$, $\tau\ne 0$, and $\Gd\ne \pm 1$. Then a point $z=\Go(\Gz)$ traces an ellipse $L_0$
whenever the point $\Gz$ traverses the unit circle $\CC_0$.

\subsection{$n=1$: a slit map}\label{s3.2}

Let $z=\Go(\Gz)$  be a conformal map that transforms the two-sided segment $l_0=l_0^+\cup l_0^-$, $l_0^\pm=[-1,1]^\pm$, into the contour $L_0$
such that the point $\Gz=\infty$ falls into the infinite point of the $z$-plane. Such a map is defined up to a real 
parameter, and it is assumed that  $\Go(\Gz)\sim c_{-1}\Gz+c_0+\ldots$, $\Gz\to\infty$, and $\I c_{-1}=0$.
Fix a single branch $q(\Gz)$ of this function in the $\Gz$-plane  cut along the segment $[-1,1]$ by the condition  $q(\Gz)\sim\Gz$,
$\Gz\to\infty$. 
Let  ${\Bbb C}_1$ and  ${\Bbb C}_2$ be two copies of the complex $\Gz$-plane with the
cut  $ l_0$. The two sheets are glued together and form a genus-0 Riemann surface  $\CR$ of the algebraic function $u^2=\Gz^2-1$ such that the loop $l_0$ 
is the symmetry 
line of the surface. Denote by $(\Gz_*,u_*)=(\bar\Gz,-u(\bar\Gz))$ the point symmetrical
to a point $(\Gz,u)$,
and let in the first sheet, $(\Gz,u)=(\Gz,q(\Gz))$.
Introduce two functions
on the surface $\CR$
\beq
\GF_1(\Gz,u)=\left\{
\begin{array}{cc}
F(\Gz), \; & (\Gz,u)\in{\Bbb C}_1,\\
\ov{F(\bar\Gz)}, \; & (\Gz,u)\in{\Bbb C}_2,\\
\end{array}
\right.
\label{3.7}
\eeq
and 
\beq
\GF_2(\Gz,u)=\left\{
\begin{array}{cc}
i\bar\tau\Go(\Gz), \; & (\Gz,u)\in{\Bbb C}_1,\\
-i\tau\ov{\Go(\bar\Gz)}, \; & (\Gz,u)\in{\Bbb C}_2.\\
\end{array}
\right.
\label{3.8}
\eeq
These two functions 
satisfy the symmetry condition
\beq
\ov{\GF_j(\Gz_*,u_*)}=\GF_j(\Gz,u), \quad (\Gz,u)\in\CR, \quad j=1,2.
\label{3.8'}
\eeq
Denote next $\GF_j^-(\Gx,v)=\ov{\GF_j^+(\Gx,v)}$,  $(\Gx,v)\in  l_0$, where  $v=u(\Gx)$,
 $\Gx\in l_0$. The two functions
are solutions to the following Riemann-Hilbert problems on the contour $l_0$:
$$
\GF_1^+(\Gx,v)-\GF_1^-(\Gx,v)=0, \quad (\Gx,v)\in  l_0,
$$
\beq
\GF_1(\Gz,u)\sim \fr{\bar\tau^\infty-\bar\tau}{\mu}c_{-1}\Gz, \quad \Gz\to\infty, \quad (\Gz,u)\in{\Bbb C}_1,
\label{3.9}
\eeq
and 
$$
\GF_2^+(\Gx,v)-\GF_2^-(\Gx,v)=2iG_0(\Gx), \quad (\Gx,v)\in l_0,
$$
\beq
\GF_2(\Gz,u)\sim i\bar\tau c_{-1}\Gz, \quad \Gz\to\infty, \quad (\Gz,u)\in{\Bbb C}_1.
\label{3.10}
\eeq
Here,  $G_0(\Gx)=\Gl_0[\R F(\Gx)-d_0']$. The  general solution of the homogeneous problem  (\ref{3.8'}), (\ref{3.9}) is given by
\beq
\GF_1(\Gz,u)=\Gb_0+\Gb_1\Gz+\Gb_2iu(\Gz), \quad (\Gz,u)\in\CR,
\label{3.11}
\eeq
where $\Gb_0$, $\Gb_1$, and $\Gb_2$ are real constants. The condition at $\infty$ is satisfied if the 
constants $\Gb_1$ and $\Gb_2$ are chosen as
\beq
\Gb_1=\fr{\tau_1^\infty-\tau_1}{\mu}c_{-1}, \quad \Gb_2=\fr{\tau_2-\tau_2^\infty}{\mu}c_{-1}.
\label{3.12}
\eeq
The solution of the inhomogeneous Riemann-Hilbert problem  (\ref{3.8'}), (\ref{3.10})
can be represented in terms of a singular integral with the Weierstrass kernel $\fr12(1+u/v)/(\Gx-\Gz)$
\beq
\GF_2(\Gz,u)=\Gg_0+\Gg_1 \Gz+i\Gg_2 u(\Gz)+\fr{1}{2\pi}\int_{l_0}G_0(\Gx)\left(1+\fr{u}{v}\right)\fr{d\Gx}{\Gx-\Gz},
\label{3.13}
\eeq
where $\Gg_0$, $\Gg_1$, and $\Gg_2$ are arbitrary real constants. The loop $ l_0$ is oriented such that the exterior of $l_0$ is on the left.
To compute the integral in formula (\ref{3.13}),  we rewrite the density $G_0(\Gz)$ as
\beq
G_0(\Gz)=\Gl_0\left[g_0(\Gx)\pm g_1(\Gx)\right], \quad \Gz=\Gx\pm i0\in l_0^\pm,
\label{3.14}
\eeq
where
\beq
g_0(\Gx)=\Gb_1\Gx, \quad g_1(\Gx)=-\Gb_2|q(\Gx)|.
\label{3.15}
\eeq
By utilizing the relation (\ref{3.14}) we deduce from formula (\ref{3.13})
\beq
\GF_2(\Gz,u)=\Gg_0+\Gg_1 \Gz+i\Gg_2 u(\Gz)+\Gl_0[\Gb_1u(\Gz)\GL_1(\Gz)-\Gb_2 \GL_2(\Gz)],
\label{3.16}
\eeq
where
\beq
\GL_1(\Gz)=\fr{1}{\pi i}\int_{-1}^1\fr{\Gx d\Gx}{\sqrt{1-\Gx^2}(\Gx-\Gz)},\quad
\GL_2(\Gz)=\fr{1}{\pi}\int_{-1}^1\fr{\sqrt{1-\Gx^2}d\Gx}{\Gx-\Gz}.
\label{3.17}
\eeq
These  integrals  are computed by the theory of residues
\beq
\GL_1(\Gz)=\fr{i\Gz}{q(\Gz)}-i,\quad 
\GL_2(\Gz)=q(\Gz)-\Gz.
\label{3.21}
\eeq
On substituting the expressions (\ref{3.21}) into (\ref{3.16}) we determine the solution to the Riemann-Hilbert problem (\ref{3.10})
explicitly
\beq
\GF_2(\Gz,u)=\Gg_0+\Gg_1\Gz+i\Gg_2u(\Gz)+\Gl_0[\Gz-q(\Gz)]\left[\Gb_2+\fr{i\Gb_1 u(\Gz)}{q(\Gz)}\right].
\label{3.22}
\eeq
Assuming that $(\Gz,u)\in{\Bbb C}_1$ we derive the formula for the conformal map
\beq
i\bar\tau\Go(\Gz)=
\Gg_0+\Gg_1\Gz+i\Gg_2q(\Gz)+\Gl_0[\Gz-q(\Gz)](i\Gb_1+\Gb_2),
\label{3.23}
\eeq
where the parameters $\Gg_1$ and $\Gg_2$ are chosen such that the function $\Go(\Gz)$ meets the condition
(\ref{2.9}) at infinity, $
\Gg_1=c_{-1}\tau_2, \quad \Gg_2=c_{-1}\tau_1. $
Finally, we substitute the parameters $\Gg_1$, $\Gg_2$, $\Gb_1$, and $\Gb_2$ 
into formula (\ref{3.23}) to deduce
\beq
\Go(\Gz)=\Gg'+c_{-1}\left(m_1\Gz+m_2\sqrt{\Gz^2-1}\right),
\label{3.25}
\eeq
where $\Gg'=-i\Gg_0/\bar\tau$ is an additive constant and
\beq
m_1=\fr{1}{\bar\tau}\left[\fr{\Gk_0}{1-\Gk_0}(\tau^\infty-\tau)-i\tau_2\right], \quad 
m_2=\fr{1}{\bar\tau}\left[-\fr{\Gk_0}{1-\Gk_0}(\tau^\infty-\tau)+\tau_1\right].
\label{3.26}
\eeq
When the points $\Gz=\Gx\pm i0$ traverses the contours $l_0^\pm$ the corresponding
points 
\beq
z=c_{-1}\left(m_1\Gx\pm im_2\sqrt{1-\Gx^2}\right)
\label{3.27}
\eeq
outline the boundary $L_0$ of a uniformly stressed inclusion $D_0$ .

\begin{figure}[t]
\centerline{
\scalebox{0.8}{\includegraphics{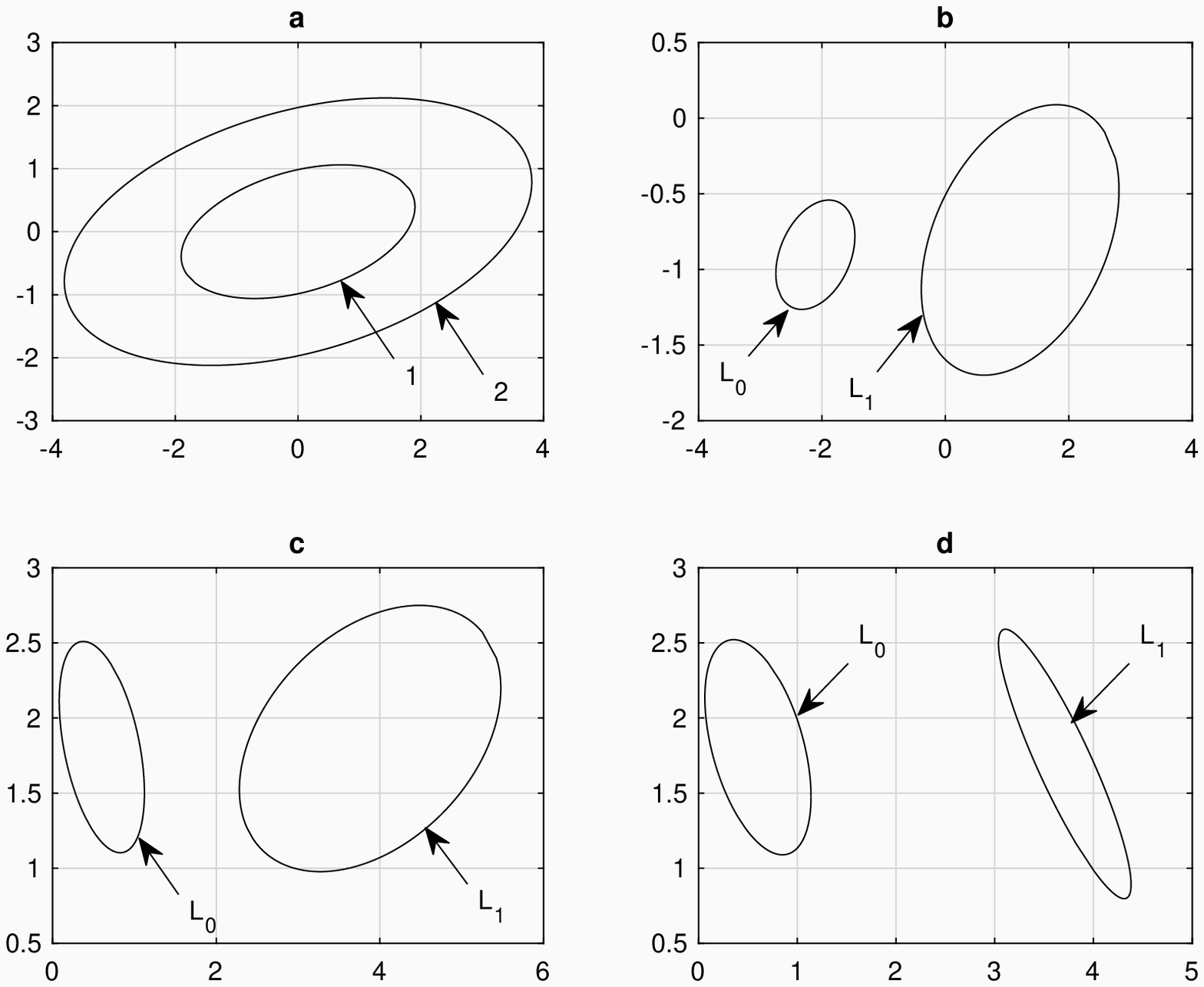}}}
\caption{Cases $n=1$ and $n=2$ when  
$\tau_1/\mu=\tau_2/\mu=-\tau_1^\infty/\mu=\tau_2^\infty/\mu=1$, $\Gg'=0$,  and $c_{-1}=1$.
(a): $n=1$. Line 1 is generated by the slit map, and line 2 is generated by the circular map.
(b) -- (d): $n=2$, $k=0.5$,  $\Gz_{\infty}=0.3k$, $a_0=\Gr_0=0$. (b): $\Gk_0=\Gk_1=5$.
(c): $\Gk_0=0.5$, $\Gk_1=5$.
(d): $\Gk_0=0.5$, $\Gk_1=0.2$.} 
\label{fig1}
\end{figure}

In Fig. \ref{fig1}a, we show the profiles of a single inclusion generated
by the slit map (line 1) and the circular map (line 2) for $\Gk=5$, $c_{-1}=1$,
and the loading parameters
$\tau_1/\mu=1$, $\tau_2/\mu=1$, $\tau_1^\infty/\mu=-1$, and $\tau_2^\infty/\mu=1$.
The contours are cofocal ellipses and coincide if $c_{-1}^s=\fr12 c_{-1}^c$, where
$c_{-1}^s$ and $c_{-1}^c$ are the scaling parameter $c_{-1}$ for the slit and circular maps, respectively.

\setcounter{equation}{0}
\section{Two inclusions}\label{s4}

 Every doubly connected domain $D^e$ may be interpreted as the image by a conformal map $z=\Go(\Gz)$
 of a slit domain $\CD^e$, the extended $\Gz$-plane cut along two segments  $l_0=[-1,-k]^+\cup[-k,-1]^-$ and 
 $l_1=[k,1]^+\cup[1,k]^-$,
 $0<k<1$. Moreover, it is always possible to choose the map such that the preimage $\Gz_\infty$
  of the infinite
  point in the $z$-plane falls into the open segment $(-k,k)$, $\Gz_\infty\in(-k.k)$. We select the single branch $q(\Gz)$ of the 
  function $u^2=(\Gz^2-1)(\Gz^2-k^2)$ in the $\Gz$-plane with the cuts
  $l_0$ and $l_1$ by the condition $q(\Gz)\sim \Gz^2$, $\Gz\to\infty$. On the banks of the cuts
   the function $q(\Gz)$ is pure imaginary,
  $q(\Gz)=\mp i|q(\Gz)|$, $\Gz\in l_0^\pm=[-1,-k]^\pm$, and  $q(\Gz)=\pm i|q(\Gz)|$, $\Gz\in l_1^\pm=[k,1]^\pm$.
It will be convenient to deal with a new function 
\beq
F_0(\Gz)=F(\Gz)-\fr{c}{\Gz-\Gz_\infty},\quad c=\fr{\bar\tau^\infty-\bar\tau}{\mu}c_{-1}.
\label{4.3}
\eeq 
 The new function has a removable singularity at the point $\Gz=\Gz_\infty$
 and satisfies the boundary condition
 \beq
 \I F_0(\Gz)=a_j-\I\left(\fr{c}{\Gz-\Gz_\infty}\right), \quad \Gz\in l_j, \quad j=0,1.
 \label{4.4}
 \eeq
Let $\CR$ be the  elliptic surface of the algebraic function $u^2=(\Gz^2-1)(\Gz^2-k^2)$ symmetric with respect to the 
contour $l=l_0\cup l_1$.
Similarly to the case $n=1$ we introduce the function $\GF_1(\Gz,u)$ on the surface $\CR$ by 
 \beq
 \GF_1(\Gz,u)=\left\{
\begin{array}{cc}
F_0(\Gz), \; & (\Gz,u)\in{\Bbb C}_1,\\
\ov{F_0(\bar\Gz)}, \; & (\Gz,u)\in{\Bbb C}_2.\\
\end{array}
\right.
\label{4.5}
\eeq 
This function is analytic everywhere in the surface $\CR$ including the two points with the affix
$\Gz=\Gz_\infty$ except for the contours $l_0$ and $l_1$  and is symmetric with respect to 
these loops.
 The limit values of the function $\GF_1(\Gz,u)$
on the cuts  $l_j$
satisfy the Riemann-Hilbert
boundary condition
\beq
\GF^+_1(\Gx,v)-\GF_1^-(\Gx,v)=2i
\left[a_j-\I\left(\fr{c}{\Gx-\Gz_\infty}\right)\right], \quad (\Gx,v)\in l_j, 
 \quad j=0,1.
\label{4.6}
\eeq
The general solution of this problem can be written in the form
\beq
\GF_1(\Gz,u)=\Gb_0+\Gb_1\Gz+\fr{1}{2\pi}\sum_{j=0}^{1}
\int_{l_j}
\left[a_j-\I\left(
\fr{c}{\Gx-\Gz_\infty}\right)\right]
\left(1+\fr{u}{v}\right)\fr{d\Gx}{\Gx-\Gz},
\label{4.7}
\eeq
where $u=u(\Gz)$, $v=v(\Gx)$, and $\Gb_j$ are real constants. The function $\GF_1(\Gz,u)$ is required to be bounded at the two  infinite points $(\infty,\infty)_1$ and  
$(\infty,\infty)_2$ of the surface $\CR$; these two points are symmetric with respect to the 
contour $l$.
Since the constants $\Gb_j$ are real and $v(\Gx)$ is pure imaginary in the loops $l_j$,
the function $\GF_1(\Gz,u)$ is bounded at the infinite points if and only if
the following conditions hold:
\beq
\Gb_1=0,
\quad
\sum_{j=0}^1\int_{l_j}\left[a_j-\I\left(\fr{c}{\Gx-\Gz_\infty}\right)\right]\fr{d\Gx}{v}=0.
 \label{4.9}
 \eeq
The symmetry of the contours $l_0$ and $l_1$ and the second condition in (\ref{4.9})
imply
\beq
a_0-a_1=\left(\int_k^1\fr{d\Gx}{|q(\Gx)|}\right)^{-1}
\left(\int_{-1}^{-k}-
\int_{k}^{1}\right)\I\left(\fr{c}{\Gx-\Gz_\infty}\right)\fr{d\Gx}{|q(\Gx)|}.
\label{4.10}
\eeq
If the constants $a_0$ and $a_1$ satisfy the condition (\ref{4.10}), then the function
\beq
\GF_1(\Gz,u)=\Gb_0+\fr{u}{2\pi}\sum_{j=0}^{1}
\int_{l_j}
\left[a_j-\I\left(
\fr{c}{\Gx-\Gz_\infty}\right)\right]
\fr{d\Gx}{(\Gx-\Gz)v}
\label{4.11}
\eeq
is bounded at infinity and defines the general solution of the problem (\ref{4.6}).  For the next step we need 
to determine the function $F_0(\Gz)$
 on the sides of the loops $l_0$ and $l_1$,
$$
 \I F_0(\Gx\pm i0)=a_j-\I\left(\fr{c}{\Gz-\Gz_\infty}\right), \quad \R F_0(\Gx\pm i0)=\Gb_0\pm (-1)^j g_1(\Gx), 
 $$
 \beq\Gx\pm i0 \in l^\pm_j, \quad j=0,1.
\label{4.12}
\eeq
Here,
\beq
g_1(\Gx)=\fr{|q(\Gx)|}{\pi}\sum_{j=0}^1 (-1)^j
\int_{l_j^+}\left[a_j-\I\left(\fr{c}{\Gn-\Gz_\infty}\right)\right]
\fr{d\Gn}{|q(\Gn)|(\Gn-\Gx)}.
\label{4.13}
\eeq
Note that the first formula in (\ref{4.12}) is consistent with the boundary condition (\ref{4.4}).

Consider now the second Schwarz problem (\ref{2.11}), (\ref{2.8}) 
for the conformal mapping $\Go(\Gz)$. As with the first Schwarz problem, we separate the singular part from the 
function $\Go(\Gz)$ by writing
\beq
\Go(\Gz)=\fr{c_{-1}}{\Gz-\Gz_\infty}+\Go_0(\Gz).
\label{4.14}
\eeq
With this splitting, the second Schwarz problem becomes
\beq
\I[i\bar\tau\Go_0(\Gz)]=G_j(\Gz)+ \Gr_j, \quad \Gz\in l_j,\quad j=0,1,
\label{4.22}
\eeq
where
$$
G_j(\Gz)=\Gl_j[g_0(\Gx)\pm(-1)^jg_1(\Gx)], \quad \Gz=\Gx\pm i0 \in l_j^\pm,
$$
\beq
\Gr_j=\Gl_j\Gr_j', \quad \Gr_j'=\Gb_0-d_j',\quad
g_0(\Gx)=\R\left[
\fr{c_{-1}}{\Gx-\Gz_\infty}\left(
\fr{\bar\tau^\infty}{\mu}-\fr{\bar\tau}{\mu_0}\right)\right].
\label{4.23}
\eeq
In order to rewrite the Schwarz problem (\ref{4.22}) as a Riemann-Hilbert problem, we introduce 
a new function on the Riemann surface $\CR$ by
\beq
\GF_2(\Gz,u)=\left\{
\begin{array}{cc}
i\bar\tau\Go_0(\Gz), \; & (\Gz,u)\in{\Bbb C}_1,\\
-i\tau\ov{\Go_0(\bar\Gz)}, \; & (\Gz,u)\in{\Bbb C}_2.\\
\end{array}
\right.
\label{4.24}
\eeq 
The function $\GF_2(\Gz,u)$ is analytic everywhere on the surface $\CR$ (including the points $\Gz_\infty$ and $\infty$) apart from the loops
 $l_j$, and on its banks it satisfies the boundary condition
 \beq
 \GF_2^+(\Gx,v)-  \GF_2^-(\Gx,v)=2i[G_j(\Gx)+\Gr_j], \quad \Gx\in l_j, \quad j=0,1.
 \label{4.25}
 \eeq
The general solution to this problem is 
 \beq
 \GF_2(\Gz,u)=\Gg_0+\Gg_1\Gz+\fr{1}{2\pi}\sum_{j=0}^{1}\int_{l_j}[G_j(\Gn)+\Gr_j]\left(1+\fr{u}{v}\right)\fr{d\Gn}{\Gn-\Gz}.
 \label{4.26}
 \eeq
The function $\GF_2(\Gz,u)$ is bounded at the two infinite
points of the surface if and only if
\beq
\Gg_1=0,\quad 
\Gr_0-\Gr_1=-\left(\int_k^1\fr{d\Gx}{|q(\Gx)|}\right)^{-1}
\left(\Gl_0\int_{-1}^{-k}-
\Gl_1\int_{k}^{1}\right)\fr{g_0(\Gx)d\Gx}{|q(\Gx)|}.
 \label{4.27}
 \eeq

The solution of the Riemann-Hilbert problem (\ref{4.25}) on the first sheet of the surface $\CR$ is
the function $i\bar\tau \Go_0(\Gz)$. By integrating over the contours $l_j^+$ and $l_j^-$ in 
 (\ref{4.26}) and employing the first relation in (\ref{4.23}) we eventually deduce
 \beq
 i\bar\tau\Go_0(\Gz)=\Gg_0+\sum_{j=0}^{1}(-1)^j\Gl_j
 \int_{l_j^+}
  \left\{g_1(\Gn)+
  \fr{i q(\Gz)}{|q(\Gn)|}
 [g_0(\Gn)+\Gr_j'] \right\}
 \fr{d\Gn}{\Gn-\Gz}.
 \label{4.30}
 \eeq
The contours $L_0$ and $L_1$ can be reconstructed by letting a point $\Gz$ 
 run along the loops $l_0$ and $l_1$. On using the Sokhotski-Plemelj  formulas 
 and the relation (\ref{4.14}) we derive
 $$
 \Go(\Gx\pm i0)=\Gg'+\fr{c_{-1}}{\Gx-\Gz_\infty}-\fr{i\mu}{\pi\bar\tau}\left\{
 \sum_{j=0}^{1}(-1)^j\tilde\Gl_j\int_{l_j^+}\left[ g_1(\Gn)\pm (-1)^m\fr{|q(\Gx)|}{|q(\Gn)|}(g_0(\Gn)+\Gr_j')\right] \fr{d\Gn}{\Gn-\Gx}
 \right.
$$
\beq
\left.
+\pi i\Gl_m[g_0(\Gx)+\Gr_m'\pm(-1)^mg_1(\Gx)]
\right\},
\quad \Gz=\Gx\pm i0\in l_m^\pm, \quad m=0,1,
\label{4.31}
\eeq 
 where $\tilde\Gl_j=\Gk_j/(1-\Gk_j)$, $\Gg'=-i\Gg_0/\bar\tau$ is an arbitrary complex constant.
 
 The conformal map, in addition to the free additive complex constant $\Gg'$, the scaling
 parameter $c_{-1}=c_{-1}'+ic_{-1}''$, and the parameters $a_0$ and $\Gr_0$, is defined up to two free real parameters, $k\in (-1,1)$ and $\Gz_\infty\in(-k,k)$.
These parameters have to be chosen such that
 the contours $L_0$  and $L_1$ outline two disjoint domains $D_0$ and $D_1$.
  Figures 1 (b -- d) sample contours $L_0$ and $L_1$ for different choices of the parameters
  $\Gk_0$ and $\Gk_1$, while the other parameters of the problem are kept the same. 
  
To reconstruct the profiles of the inclusions, one needs to compute some integrals.  The integrals in (\ref{4.10})
and (\ref{4.27}) are calculated by the Gauss quadrature formulas
\beq
\int_a^b\fr{h(\Gx)d\Gx}{\sqrt{(\Gx-a)(b-\Gx)}}\approx\fr{\pi}{N}\sum_{j=1}^Nh(\Gx_j),
\label{4.32} 
 \eeq
 where
 \beq
 \Gx_j=\Gd_++\Gd_-x_j, \quad \Gd_\pm=\fr{b\pm a}{2}, \quad x_j=\cos\fr{(2j-1)\pi}{2N}.
 \label{4.33}
 \eeq
 The evaluation of the function $g_1(\Gx)$ and $\Go(\Gx\pm i0)$ ($\Gx\pm i0\in l_m^\pm$)  by formulas 
 (\ref{4.13}) and (\ref{4.31}), respectively, 
 requires computing singular integrals with the Cauchy kernel of the form
\beq
S(\Gx)=\int_a^b\fr{h(\Gn)d\Gn}{\sqrt{(\Gn-a)(b-\Gn)}(\Gn-\Gx)}, \quad a< \Gx< b.
\label{4.34} 
 \eeq 
On making the substitutions
\beq
\Gx=\Gd_++\Gd_-x, \quad  \Gn=\Gd_++\Gd_-y,\quad h(\Gn)=\tilde h(y), \quad -1<y<1,\quad -1<x<1,
\label{4.35}
\eeq
we expand the function $\tilde h(y)$ in terms of the Chebyshev polynomials of the first kind
\beq
\tilde h(y)=\sum_{m=0}^\infty \Ga_m T_m(y).
\label{4.36}
\eeq 
 Here,
 \beq
 \Ga_0=\fr{1}{\pi}\int_{-1}^1\fr{\tilde h(y)dy}{\sqrt{1-y^2}},
 \quad \Ga_m=\fr{2}{\pi}\int_{-1}^1\fr{\tilde h(y)T_m(y)dy}{\sqrt{1-y^2}},\quad m=1,2,\ldots.
\label{4.37}
 \eeq
 Approximately, by the Gauss formula,
 \beq
 \Ga_m\approx \fr{2}{N}\sum_{j=1}^N \tilde h(x_j)\cos\fr{(2j-1)\pi m}{2N}, \quad m=1,2,\ldots.
 \label{4.38}
 \eeq
 By substituting the expansion (\ref{4.36}) into formula (\ref{4.34}) and employing the relation
 \beq
 \int_{-1}^1\fr{T_m(y)dy}{\sqrt{1-y^2}(y-x)}=\left\{
 \begin{array}{cc}
 0, & m=0,\\
 \pi U_{m-1}(x), & m=1,2,\ldots,\\
 \end{array}
 \right.\quad -1<x<1,
 \label{4.39}
 \eeq
 where $U_m(x)$ is the Chebyshev polynomial of the second kind,
 we finally obtain the formula for the singular integral (\ref{4.34})
 \beq
 S(\Gx)=\fr{\pi}{\Gd_-}\sum_{m=1}^\infty \Ga_m U_{m-1}\left(\fr{\Gx-\Gd_+}{\Gd_-}\right).
 \label{4.40}
 \eeq

To reconstruct two uniformly stressed inclusions symmetric with respect to the origin,
we choose $\Gk_0=\Gk_1$ and $\Gz_\infty=0$. Then
  \beq
  a_1=-a_0=c''_{-1}\left(\int_k^1\fr{d\Gx}{|q(\Gx)|}\right)^{-1}\int_k^1\fr{d\Gx}{\Gx|q(\Gx)|}.
  \label{4.41}
  \eeq
  The functions $g_1(\Gx)$ and $g_0(\Gx)$ given by (\ref{4.13}) and (\ref{4.23})
  are both odd and have the form
  $$
  g_0(\Gx)=\fr{c_*}{\Gx}, \quad c_*=\R\left[c_{-1}\left(\fr{\bar\tau^\infty}{\mu}-\fr{\bar\tau}{\mu_0}\right)\right],
  $$
  \beq
  g_1(\Gx)=\fr{2\Gx|q(\Gx)|}{\pi}\int_k^1\left(\fr{c'}{\Gn}-a_1\right)\fr{d\Gn}{|q(\Gn)|(\Gn^2-\Gx^2)}.
  \label{4.42}
  \eeq
This simplifies formula (\ref{4.27}), and the constants $\Gr_0$ and $\Gr_1$ have the form
 \beq
  \Gr_1=-\Gr_0=-\fr{\mu_0 c_*}{1-\Gk_0}\left(\int_k^1\fr{d\Gx}{|q(\Gx)|}\right)^{-1}\int_k^1\fr{d\Gx}{\Gx|q(\Gx)|}.
  \label{4.43}
  \eeq
 Our next step is to determine the conformal map  function $z=\Go(\Gz)$ on the boundary of the loops 
$l_0$ and $l_1$. Without loss we may drop the additive constant $\Gg_0$. Since $g_j(\Gx)=-g_j(-\Gx)$,
$j=0,1$, on using the representation  (\ref{4.14}) and the Sokhotski-Plemelj formulas 
we deduce 
 $$
\Go(-\Gx\mp i0)=-\fr{c_{-1}}{\Gx}+\fr{i\mu\tilde\Gl_0}{\bar\tau}\left\{
 \mp i g_1(\Gx)+ig_0(\Gx)+i\Gr_1'
\right.
$$ 
\beq
\left.
+\fr{1}{\pi}\sum_{j=0}^1(-1)^j\int_{k}^{1}
  \left[-g_1(\Gn)\pm\fr{|q(\Gx)|}{|q(\Gn)|}
 (g_0(\Gn)+\Gr_1') \right]
 \fr{d\Gn}{\Gn-(-1)^j\Gx}
 \right\}, \quad -\Gx\mp i0\in l_0^\mp.
 \label{4.45}
 \eeq 
Similarly, we let $\Gz=\Gx\pm i0\in l_1^\pm$ and derive
that
$
\Go(\Gx\pm i0)=-\Go(-\Gx\mp i0),\quad \Gx\pm i0\in l_1^\pm.
$
This implies that the two inclusions are symmetric with respect to the origin.

\begin{figure}[t]
\centerline{
\scalebox{0.8}{\includegraphics{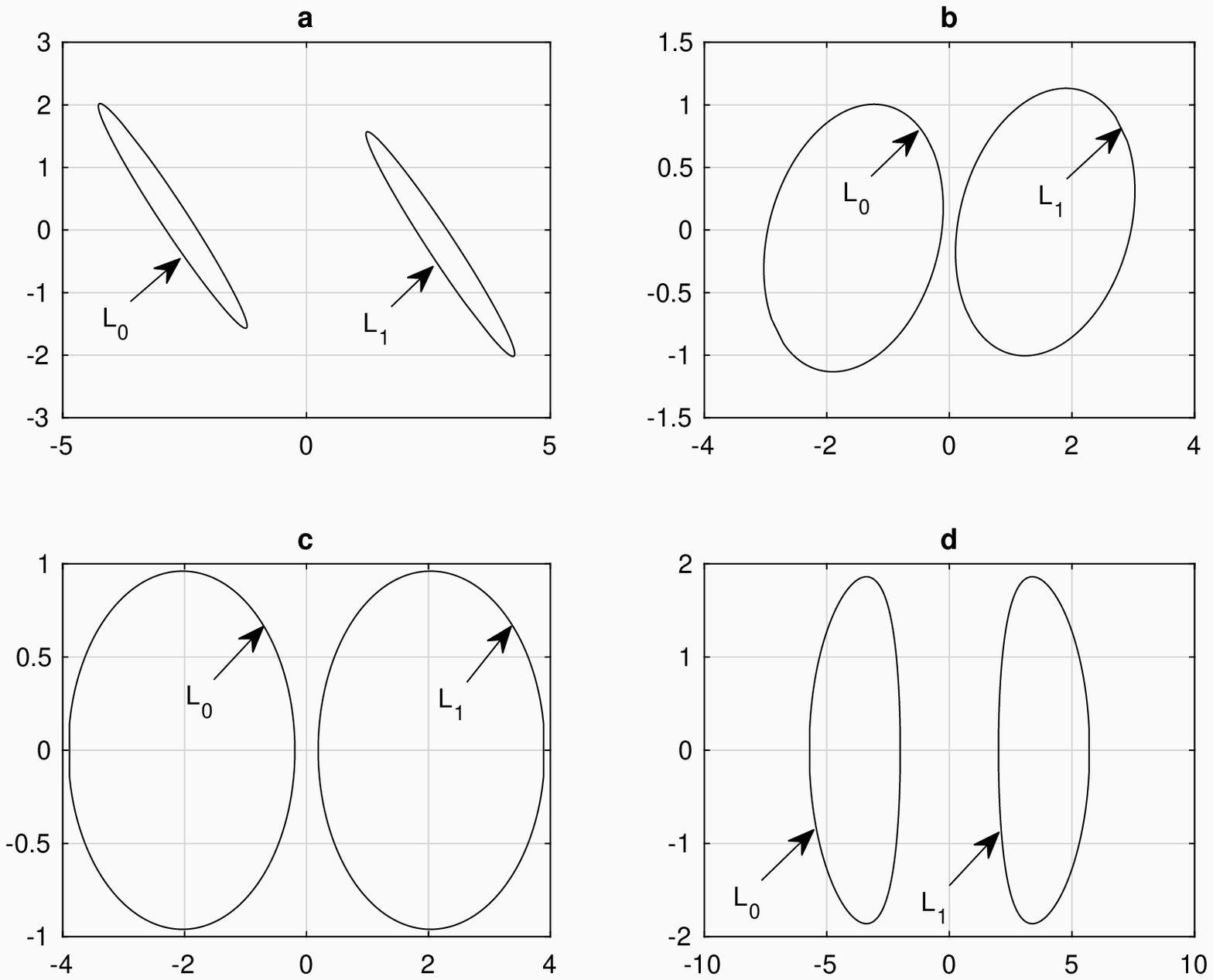}}}
\caption{Two symmetric inclusions ($\Gz_\infty=0$, $a_0=-a_1$, $\Gr_0=-\Gr_1$) 
when $\tau_1/\mu=-\tau_1^\infty/\mu=1$
and $c_{-1}=1$.
(a): $\tau_2/\mu=\tau_2^\infty/\mu=1$, $k=0.2$, and $\Gk_0=\Gk_1=0.1$.
(b): $\tau_2/\mu=\tau_2^\infty/\mu=1$, $k=0.5$, and $\Gk_0=\Gk_1=2$.
(c) -- (d): $\tau_2=\tau_2^\infty=0$, $k=0.35$, and $\Gk_0=\Gk_1=1000$. (c): the solvability conditions (\ref{4.41}) and
(\ref{4.43}) are satisfied. 
(d):  $a_0=-a_1=0$ and $\Gr_0=-\Gr_1=0$ -- the solvability conditions (\ref{4.41}) and
(\ref{4.43}) are not satisfied.
} 
\label{fig2}
\end{figure}

In Figures \ref{fig2} (a -- d) we outline the profiles  of two symmetric inclusions.   Equations  (\ref{4.41}) and (\ref{4.43}) 
fix the parameters $a_j=\Gk_j d_j''-b_j$ and $\Gr_j=\Gl_j(\Gb_0-d_j')$.  These equations are the solvability conditions for 
the two Schwarz problems (\ref{2.10}) and (\ref{2.11}) in the doubly-connected slit domain or, equivalently,   
the solvability conditions for  the inverse problem of elasticity (\ref{2.4}).
If these conditions are violated, then the resulting conformal map generates different inclusions' profiles.
In Figures 2 (c), we show the contours $L_0$ and $L_1$ when the parameters  $a_0=-a_1$ and $\Gr_0=-\Gr_1$  are determined
by  equations (\ref{4.41}) and
(\ref{4.43}), while in Figure  2(d), for the same set of the problem parameters, we outline the contours $L_0$ and $L_1$ when  $a_0=-a_1=0$ and $\Gr_0=-\Gr_1=0$, 
and the solvability conditions are violated.

\setcounter{equation}{0}

\setcounter{equation}{0}

\section{Three inclusions}

 If $D^e$ is any triply connected domain, then there exists a conformal map $z=\Go(\Gz)$
  such that it transforms the exterior $\CD^e$ of three two-sided segments 
$l_0=[-1,-k]^+\cup[-k,-1]^-$,
$l_1=[k_1,k_2]^+\cup[k_2,k_1]^-$, and $l_2=[k,1]^+\cup[1,k]^-$ into the domain $D^e$. Here,
$0<k<1$ and $-k<k_1<k_2<k$. 
The point $z=\infty$ is the image of a certain point $\Gz_\infty=\Gz_\infty'+i\Gz_\infty''$ of the parametric $\Gz$-plane, and, in general,
the parameters $\Gz_\infty'$, $\Gz_\infty''$, $k$, $k_1$, and $k_2$
 cannot be prescribed and should be recovered from some additional conditions. 
 
Denote 
\beq
p(\Gz)=(\Gz^2-1)(\Gz^2-k^2)(\Gz-k_1)(\Gz-k_2)
\label{5.1}
\eeq
and fix a single branch $q(\Gz)$ of the function $u^2=p(\Gz)$ in the $\Gz$-plane with the cuts
$l_0$, $l_1$, and $l_2$ by the condition $q(\Gz)\sim \Gz^3$, $\Gz\to\infty$. This branch is pure imaginary
on the cuts' sides, $q(\Gz)=\pm(-1)^ji|q(\Gz)|$, $\Gz\in l_j^\pm$.
As before, we aim to reduce the Schwarz problems (\ref{2.10}) and (\ref{2.11})
to two Riemann-Hilbert problems on the Riemann surface $\CR$  of the algebraic function $u^2=p(\Gz)$.
In the case of three inclusions,
$\CR$  is a hyperelliptic genus-2 surface symmetric with respect to the contour $l=l_0\cup l_1\cup l_2$.

\subsection{$\Gz_\infty$ is a finite point}

Suppose first that the preimage of the point $z=\infty$ is a finite point $\Gz_\infty=\Gz_\infty'+i\Gz_\infty''$. 
Then the conformal mapping $\Go(\Gz)$ has a simple pole at the point $\Gz_\infty$, and its expansion
at this point has the form (\ref{2.8}).
The general solution of the 
corresponding Riemann-Hilbert problem (\ref{4.6}) for the function $\GF_1(\Gz,u)$
introduced in (\ref{4.5}) is given by 
\beq
\GF_1(\Gz,u)=\Gb_0+\Gb_1\Gz+\Gb_2\Gz^2+\fr{1}{2\pi}\sum_{j=0}^{2}
\int_{l_j}
\left[a_j-\I\left(
\fr{c}{\Gx-\Gz_\infty}\right)\right]
\left(1+\fr{u}{v}\right)\fr{d\Gx}{\Gx-\Gz},
\label{5.2}
\eeq
where $\Gb_j$ are real constants. 
In general, this solution is not bounded at infinity.
The necessary and sufficient conditions for the solution to be bounded at the infinite points of the surface $\CR$ are
\beq
\Gb_1=\Gb_2=0,\quad
 \sum_{j=0}^{2}\int_{l_j}\left[a_j-\I\left(\fr{c}{\Gx-\Gz_\infty}\right)\right]\fr{\Gx^{m-1}d\Gx}{v}=0, \quad m=1,2.
 \label{5.3}
 \eeq
Denote the integrals
 \beq
 I_{mj}=\int_{l_j^+}\fr{\Gx^{m-1}d\Gx}{|q(\Gx)|}, \quad j=0,1,2, \quad m=1,2.
 \label{5.4}
 \eeq
In terms of these integrals and the constant $a_0$ the other constants $a_1$ and $a_2$ are expressed 
as
\beq
a_1=\fr{J_2I_{12}-J_1I_{22}}{\GD},\quad
a_2=\fr{J_2I_{11}-J_1I_{21}}{\GD}.
 \label{5.5}
 \eeq
 Here,
 $$
 J_m =J_{m0}-J_{m1}+J_{m2}-a_0I_{m0}, \quad \GD=I_{11}I_{22}-I_{12}I_{21},
 $$
 \beq
 J_{mj}=\int_{l_j^+}\I\left(\fr{c}{\Gx-\Gz_\infty}\right)\fr{\Gx^{m-1}d\Gx}{|q(\Gx)|}.
 \label{5.6}
 \eeq
By utilizing formulas (\ref{4.3}), (\ref{4.5}), (\ref{5.2}), and (\ref{5.3}) we derive
\beq
F(\Gz)=\Gb_0+\fr{c}{\Gz-\Gz_\infty}+\fr{q(\Gz)}{2\pi}\sum_{j=0}^{2}
\int_{l_j}
\left[a_j-\I\left(
\fr{c}{\Gx-\Gz_\infty}\right)\right]
\fr{d\Gx}{(\Gx-\Gz)v}.
\label{5.7}
\eeq
As in the case $n=2$ we employ the function $\GF_2(\Gz,u)$ introduced in (\ref{4.24}).
The function  $\GF_2(\Gz,u)$ solves the Riemann-Hilbert problem on the genus-2 Riemann surface
\beq
 \GF_2^+(\Gx,v)-  \GF_2^-(\Gx,v)=2i[G_j(\Gx)+\Gr_j], \quad \Gx\in l_j, \quad j=0,1,2.
 \label{5.8}
 \eeq
where $\Gr_j=\Gl_j\Gr_j'$, $\Gr_j'=\Gb_0-d_j'$, and $G_j(\Gz)$ is given by (\ref{4.23}) and 
\beq
g_1(\Gx)=\fr{|q(\Gx)|}{\pi}\sum_{j=0}^2(-1)^j
\int_{l_j^+}\left[a_j-\I\left(\fr{c}{\Gn-\Gz_\infty}\right)\right]
\fr{d\Gn}{|q(\Gn)|(\Gn-\Gx)}.
\label{5.9}
\eeq
The solution to the problem (\ref{5.8}) 
can be represented in the first sheet ${\Bbb C}_1$ as
\beq
 i\bar \tau\Go_0(\Gz)=\Gg_0+\fr{1}{\pi}\sum_{j=0}^{2}(-1)^j\Gl_j
 \int_{l_j^+}
  \left\{g_1(\Gn)-
  \fr{i q(\Gz)}{|q(\Gn)|}
 [g_0(\Gn)+\Gr_j'] \right\}
 \fr{d\Gn}{\Gn-\Gz}.
 \label{5.10}
 \eeq
This function is bounded at infinity if and only if the constants $\Gr_j$ are selected to be
$$
\Gr_1=\fr{K_1I_{22}-K_2I_{12}}{\GD},\quad
\Gr_2=\fr{K_1I_{21}-K_2I_{11}}{\GD},
 $$
 \beq
 K_m =K_{m0}-K_{m1}+K_{m2}+\Gr_0 I_{m0}, \quad 
 K_{mj}=\Gl_j\int_{l_j^+}\fr{g_0(\Gx)\Gx^{m-1}d\Gx}{|q(\Gx)|},
 \label{5.12}
 \eeq
and $\Gr_0$ is a free parameter.

The inclusions' profiles $L_0$, $L_0$, and $L_2$ are described by the 
function $\Go(\Gz)$ when a point $\Gz$ traverses the loops $l_0$, $l_1$, and $l_2$, respectively.
On these contours, the function $\Go(\Gz)$ 
is determined by
 $$
 \Go(\Gx\pm i0)=\Gg'+\fr{c_{-1}}{\Gx-\Gz_\infty}-\fr{i\mu}{\pi\bar\tau}\left\{
 \sum_{j=0}^{2}(-1)^j\tilde\Gl_j\int_{l_j^+}\left[ g_1(\Gn)\pm (-1)^m\fr{|q(\Gx)|}{|q(\Gn)|}(g_0(\Gn)+\Gr_j')\right] 
\right.
$$
\beq
\left.
\times\fr{d\Gn}{\Gn-\Gx}+\pi i\tilde\Gl_m[g_0(\Gx)+\Gr_m'\pm(-1)^mg_1(\Gx)]
\right\},
\quad \Gz=\Gx\pm i0\in l_m^\pm, \quad m=0,1,2.
\label{5.13}
\eeq 

\begin{figure}[t]
\centerline{
\scalebox{0.8}{\includegraphics{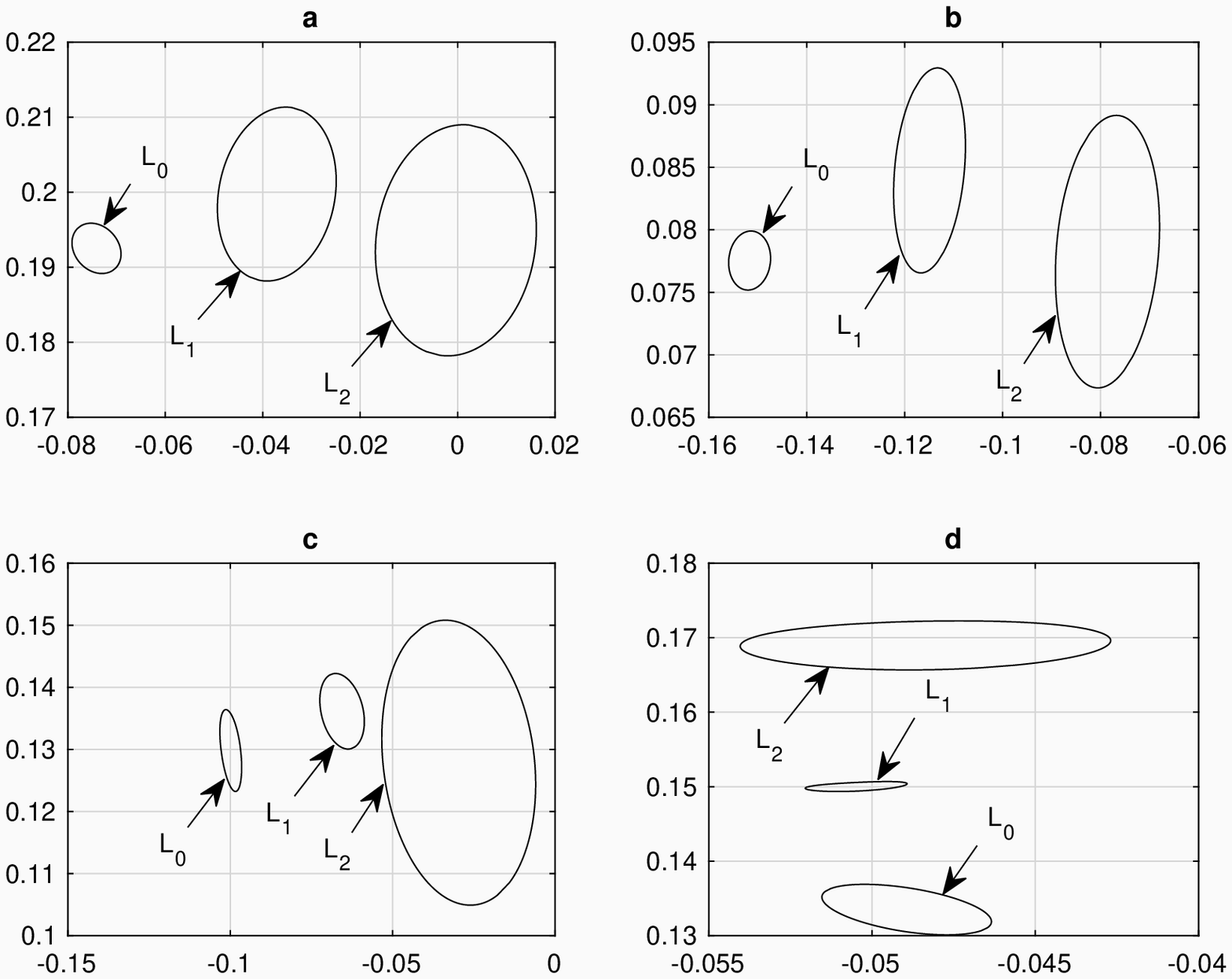}}}
\caption{Three inclusions when  $\Gz_\infty$ is a finite point, 
$\tau_1/\mu=-\tau_1^\infty/\mu=1$, $k_1=-0.1$, $a_0=0$, $\Gr_0=0$, $\Gg'=0$, 
and $c_{-1}=1$. 
(a): $\tau_2=\tau_2^\infty=0$, $\Gk=0.1$, $k=0.5$,  $k_2=0.1$, and $\Gz_\infty=5i$.
(b): $\tau_2=\tau_2^\infty=0$, $\Gk=0.2$, $k=0.5$, $k_2=0.1$, and $\Gz_\infty=5i$.
(c): $\tau_2=\tau_2^\infty=1$, $\Gk=0.3$, $k=0.6$,   $k_2=0.1$, and $\Gz_\infty=1+2i$.
(d): $\tau_2=\tau_2^\infty=1$, $\Gk=0.1$, $k=0.5$,   $k_2=0.2$, and $\Gz_\infty=5+5i$.}  
\label{fig3}
\end{figure} 
The conformal mapping $z=\Go(\Gz)$ described by (\ref{5.13}) in addition to the scaling parameter
$c_{-1}=c_{-1}'+ic_{-1}''$  has seven real free parameters $\Gz_{\infty}'$, $\Gz_{\infty}''$,
$k$, $k_1$, $k_2$, $a_0$, and $\Gr_0$.
Sample profiles of three uniformly stressed inclusions are shown in Figures 3 (a)--(d) in the
case when the preimage of the infinite point  is  a finite point $\Gz_\infty$ of the parametric
$\Gz$-plane.

\subsection{$\Gz_\infty$ is the infinite point}

If the preimage of the infinite point $z=\infty$ is the infinite point $\Gz=\infty$, 
then in a neighborhood of the infinite point
the functions $\Go(\Gz)$ and $F(\Gz)$ have expansions (\ref{2.9}) and (\ref{2.13}),
respectively. 
Instead of the function (\ref{4.3}) we introduce a new function $F_0(\Gz)$ by
\beq
F_0(\Gz)=F(\Gz)-c\Gz, \quad c=c'+ic''.
\label{5.14}
\eeq
Then the function $F_0(\Gz)$ has a removable singularity at the infinite point. My making use of this function
similarly to the case when $\Gz_\infty$ is a finite point we derive the 
 solution to the first Schwarz problem (\ref{2.10}), (\ref{2.13}). It is
\beq
F(\Gz)=\Gb_0+c\Gz+\fr{q(\Gz)}{2\pi}\sum_{j=0}^{2}
\int_{l_j}
\fr{(a_j-c''\Gx)d\Gx}{(\Gx-\Gz)v}.
\label{5.15}
\eeq
Only one constant say, $a_0$, can be chosen arbitrarily. The solution $F(\Gz)$
is bounded at infinity if and only if the other two constants are chosen as
\beq
a_1=\fr{\hat J_2I_{12}-\hat J_1I_{22}}{\GD}, \quad a_2=\fr{\hat J_2I_{11}-\hat J_1I_{21}}{\GD}, 
\label{5.16}
\eeq
where
\beq
\hat J_m=c''\sum_{j=0}^2(-1)^jI_{m+1 j}-a_0I_{m0},  \quad I_{sj}=\int_{l_j^+}\fr{\Gx^{s-1}d\Gx}{|q(\Gx)|}.
\label{5.17}
\eeq
In the same fashion instead of the function (\ref{4.14}) we introduce the function
\beq
\Go_0(\Gz)=\Go(\Gz)-c_{-1}\Gz
\label{5.18}
\eeq
and reduce the second Schwarz problem (\ref{2.9}), (\ref{2.11}) to the following:
\beq
\I[i\bar\tau\Go_0(\Gz)]=\Gl_j[g_0(\Gx)\pm(-1)^jg_1(\Gx)+\Gr_j'], \quad \Gz=\Gx\pm i0\in l_j^\pm,
\label{5.19}
\eeq
where
$$
g_0(\Gx)=c_{*j}\Gx, \quad c_{*j}=\R\left[c_{-1}\left(\fr{\bar \tau^\infty}{\mu}-\fr{\bar\tau}{\mu_j}\right)\right],
$$
\beq
g_1(\Gx)=\fr{|q(\Gx)|}{\pi}\sum_{j=0}^2 (-1)^j
\int_{l_j^+}\fr{(a_j-c''\Gn)d\Gn}{|q(\Gn)|(\Gn-\Gx)}.
\label{5.20}
\eeq
The solution of the Schwarz problem (\ref{5.19}) is given by formula (\ref{5.10}) in the domain $\CD^e={\Bbb C}\setminus l$ and  formula (\ref{5.13}) in the contours $l_0$, $l_1$, and $l_2$.
In these formulas the functions $g_0(\Gn)$ and $g_1(\Gn)$ need to be replaced by  the ones in (\ref{5.20}),
while the constants $\Gr_j'$ are
$$
\Gr_j'=\fr{\Gr_j}{\Gl_j},\quad \Gr_1=\fr{\hat K_1 I_{22}-\hat K_2I_{12}}{\GD}, \quad \Gr_2=\fr{\hat K_1I_{21}-\hat K_2I_{11}}{\GD}, 
$$
\beq
\hat K_m=\sum_{j=0}^2(-1)^jc_{*j}\Gl_jI_{m+1 j}+\Gr_0I_{m0}, \quad m=1,2.
\label{5.21}
\eeq

To reconstruct symmetrically located uniformly stressed inclusions in the case $\Gz_\infty=\infty$,
we take $\tau_2=\tau^\infty_2$, $a_0=-a_2$,  $\Gr_0=-\Gr_2$, $k_1=-k_2$, and $\Gk_0=\Gk_1=\Gk_2$.
For $\Gz_\infty=\infty$, samples of two symmetric and nonsymmetric  inclusions are shown in Figures 4 (a), (b) and $4 (c), (d)$, respectively. 
In Figure 4 (d), the contours intersect each other, and  the set of parameters is inacceptable.

\begin{figure}[t]
\centerline{
\scalebox{0.8}{\includegraphics{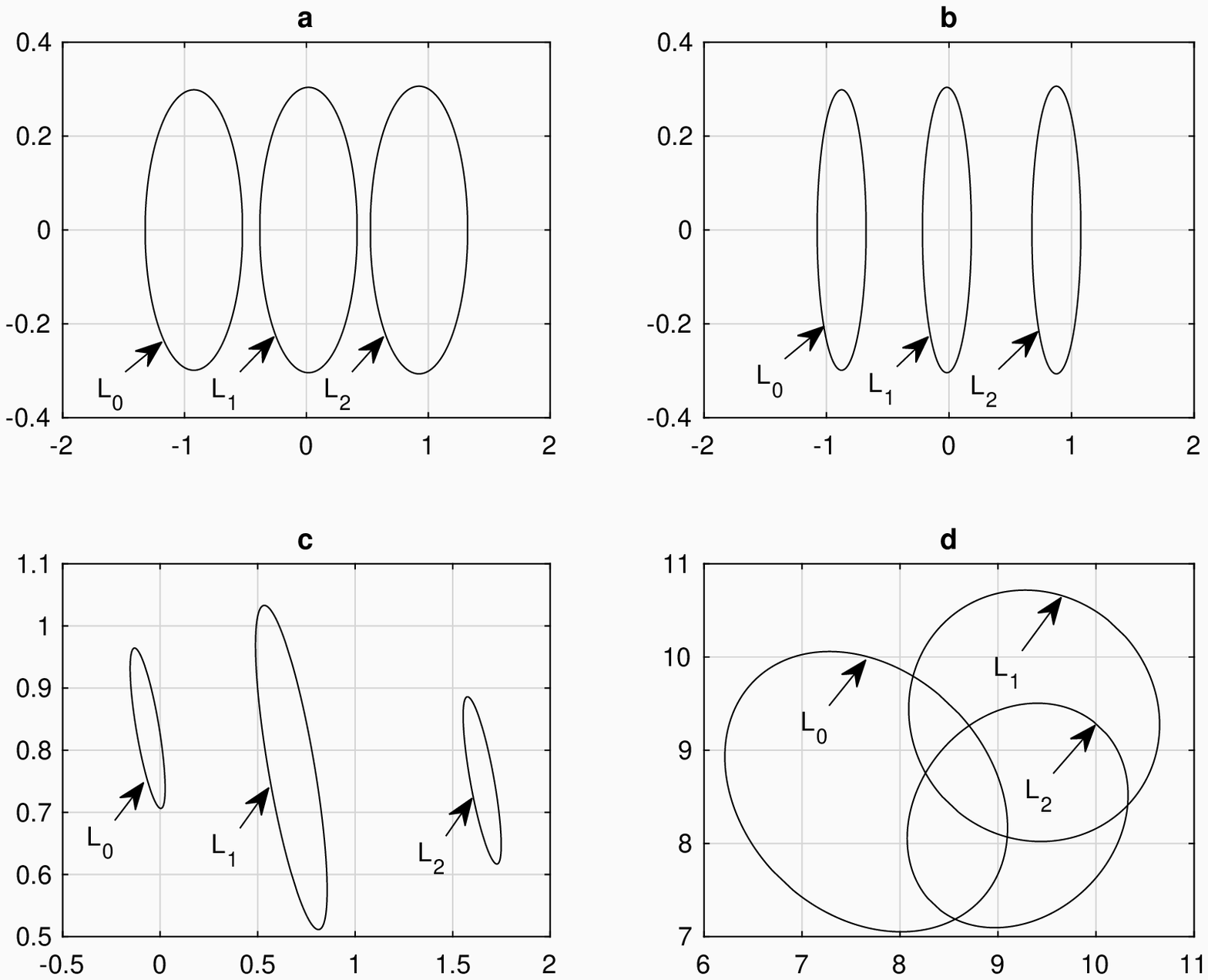}}}
\caption{Three inclusions when  $\Gz_\infty=\infty$, $\tau_1/\mu=-\tau_1^\infty/\mu=1$, $k=0.8$, $k_2=0.1$, $\Gg'=0$, 
and $c_{-1}=1$.
(a): $\tau_2=\tau_2^\infty=0$, $\Gk_0=\Gk_1=\Gk_2=2$,   $k_1=-0.1$, $a_0=-a_2$, and $\Gr_0=-\Gr_2$.
(b): $\tau_2=\tau_2^\infty=0$, $\Gk_0=\Gk_1=\Gk_2=0.5$,  $k_1=-0.1$, $a_0=-a_2$, and $\Gr_0=-\Gr_2$.
(c): $\tau_2/\mu=\tau_2^\infty/\mu=1$, $\Gk_0=\Gk_1=\Gk_2=0.3$,   $k_1=-0.3$, $a_0=0$, and $\Gr_0=0$.
(d): $\tau_2/\mu=\tau_2^\infty/\mu=1$, $\Gk_0=\Gk_1=\Gk_2=0.9$,  $k_1=-0.1$, $a_0=0$, and $\Gr_0=0$ -- no solution
for this set of parameters.}
\label{fig4}
\end{figure}

\setcounter{equation}{0}
\section{$n\ge 4$ inclusions}\label{s4}

Assume that $n\ge 4$ and $D^e={\Bbb R}^2\setminus D$ is an $n$-connected domain  ($D=\cup_{j=0}^{n-1}D_j$) such that there exists a conformal map  $z=\Go(\Gz)$
that transforms a slit domain $\CD^e$ into the domain $D^e$, and  $\CD^e$ is the exterior of $n$ slits $l_j=[k_{2j}, k_{2j+1}]^+\cup[k_{2j+1}, k_{2j}]^-$ 
($j=0,1,\ldots,n-1$) lying in the real $\Gz$-axis,   $k_0<k_1<\ldots<k_{2n-2}<k_{2n-1}$, and
$k_0=-1$, $k_{2n-1}=1$. 
When $n\ge 4$  not for  each $n$-connected domain $D^e$ there exists a conformal map that transforms the exterior 
$\CD^e$ of $n$ slits   into the domain $D^e$ such that all the slits lie in the same line. In what follows  we confine ourselves to the family of domains $D^e$ when such a map exists. Denote that family $D^e_*$. 
Notice that every doubly and triply connected domain belongs to the family $D_*^e$. 
We confine ourselves to the case when the preimage of the point $z=\infty$
is a finite point $\Gz_\infty$.

First we fix the branch $q(\Gz)$
of the   function 
\beq
p^{1/2}(\Gz)=\prod_{j=0}^{2n-1} (\Gz-k_j)^{1/2}
\label{6.1}
\eeq
in the domain $\CD^e$ by the condition $q(\Gz)\sim\Gz^n$, $\Gz\to\infty$.
On the sides of the cuts this branch is pure imaginary and
\beq
q(\Gz)|_{\Gz\in l_{n-j-1}^\pm}=\pm i(-1)^j | q(\Gz)|, \quad j=0,1,\ldots,n-1.
\label{6.2}
\eeq
As in the cases $n=2$ and $n=3$ we introduce the functions $F_0(\Gz)$ and $\Go_0(\Gz)$
by formulas (\ref{4.3}) and (\ref{4.14}), respectively. In terms of the function $\GF_1(\Gz,u)$
given by (\ref{4.5}) on the genus-$(n-1)$ Riemann surface $\CR$ of the algebraic function $u^2=p(\Gz)$
the first Schwarz problem (\ref{2.10}), (\ref{2.12}) can be equivalently written as
the Riemann-Hilbert problem
\beq
\GF^+_1(\Gx,v)-\GF_1^-(\Gx,v)=2i
\left[a_j-\I\left(\fr{c}{\Gx-\Gz_\infty}\right)\right], \quad (\Gx,v)\in l_j, 
 \quad j=0,1, \ldots, n-1.
\label{6.3}
\eeq
The general solution of this problem is 
\beq
\GF_1(\Gz,u)=\Gb_0+\fr{1}{2\pi}\sum_{j=0}^{n-1}
\int_{l_j}
\left[a_j-\I\left(
\fr{c}{\Gx-\Gz_\infty}\right)\right]
\left(1+\fr{u}{v}\right)\fr{d\Gx}{\Gx-\Gz}.
\label{6.4}
\eeq
For the function $\GF_1(\Gz,u)$  to be bounded at the two  infinite points $(\infty,\infty)_1$ and  
$(\infty,\infty)_2$ of the surface $\CR$ it is necessary and sufficient that
 \beq
 \sum_{j=0}^{n-1}\int_{l_j}\left[a_j-\I\left(\fr{c}{\Gx-\Gz_\infty}\right)\right]\fr{\Gx^{m-1}d\Gx}{v}=0, \quad m=1,2,\ldots, n-1.
 \label{6.5}
 \eeq
 The integrals 
 \beq
 A_{mj}=\int_{l_j}\fr{\Gx^{m-1}d\Gx}{v}, \quad j=1,\ldots,n-1, \quad m=1,\ldots,n-1,
 \label{6.6}
 \eeq
 form the matrix of the $A$-periods of the abelian integrals (Springer, 1956)
 \beq
 \int_{(-1,0)}^{(\Gx,u)}\fr{\Gx^{m-1}d\Gx}{v}, \quad m=1,\ldots,n-1.
 \label{6.7}
 \eeq
Select $a_0$ as a free parameter. Then the other constants $a_j$ are recovered from the system 
 of linear algebraic equations
 \beq
 \sum_{j=1}^{n-1} A_{mj} a_j=b_m, \quad m=1,\ldots,n-1.
 \label{6.8}
 \eeq
 Here,
 \beq
 b_m =
 \sum_{j=0}^{n-1}\int_{l_j}\I\left(\fr{c}{\Gx-\Gz_\infty}\right)\fr{\Gx^{m-1}d\Gx}{v}-A_{m0}a_0.
 \label{6.9}
 \eeq
Since the matrix of  the $A$-periods of basis abelian integrals is not singular, the solution to the system (\ref{4.12})
exists and is unique.
The function $v$ is pure imaginary on the sides $l_j^\pm$ of the loops $l_j$. Therefore
 employing the relations (\ref{6.2}) we represent the system (\ref{6.8}) as
 \beq
\sum_{j=1}^{n-1} (-1)^{j}I_{mj}a_j =J_m, \quad m=1,2,\ldots,n-1,
\label{6.10}
\eeq
where
 \beq
I_{mj} =\int_{l_j^+}\fr{\Gx^{m-1}d\Gx}{|q(\Gx)|}>0,
\quad J_m= 
  \sum_{j=0}^{n-1}(-1)^{j}\int_{l_j^+}\I\left(\fr{c}{\Gx-\Gz_\infty}\right)\fr{\Gx^{m-1}d\Gx}{|q(\Gx)|}-I_{m0}a_0.
 \label{6.11}
 \eeq
Since the coefficients $I_{mj}$ and  $J_m$ are real, 
the constants $a_j$ ($j=1,\ldots,n-1$) are also real.  When $(\Gz,u)\in {\Bbb C}_1$, $\GF_1(\Gz,u)=F_0(\Gz)$,
and therefore the function $F(\Gz)$ admits the following representation by quadratures:
\beq
F(\Gz)=\Gb_0+\fr{c}{\Gz-\Gz_\infty}+\fr{q(\Gz)}{2\pi}\sum_{j=0}^{n-1}\int_{l_j}\left[a_j-\I\left(\fr{c}{\Gn-\Gz_\infty}\right)\right]
\fr{d\Gn}{(\Gn-\Gz)v}.
\label{6.12}
\eeq

Proceed now to the second Schwarz problem (\ref{2.11}), (\ref{2.8}) 
equivalent to 
\beq
\I[i\bar\tau\Go_0(\Gz)]=G_j(\Gz)+ \Gr_j, \quad \Gz\in l_j,\quad j=0,1,\ldots,n-1,
\label{6.13}
\eeq
where  $\Gr_j=\Gl_j\Gr_j'$, $\Gr_j'=\Gb_0-d_j'$,
the functions $G_j(\Gz)$ and $g_0(\Gx)$ are defined in (\ref{4.23}), and 
\beq
g_1(\Gx)=\fr{|q(\Gx)|}{\pi}\sum_{j=0}^{n-1} (-1)^j
\int_{l_j^+}\left[a_j-\I\left(\fr{c}{\Gn-\Gz_\infty}\right)\right]
\fr{d\Gn}{|q(\Gn)|(\Gn-\Gx)}.
\label{6.15}
\eeq
Similarly to the cases $n=2$ and $n=3$ we rewrite the Schwarz problem (\ref{6.13}) as a Riemann-Hilbert problem
on the surface $\CR$. Ultimately, we derive 
\beq
 i\bar\tau\Go_0(\Gz)=\Gg_0+\fr{1}{\pi}\sum_{j=0}^{n-1}(-1)^j\Gl_j
 \int_{l_j^+}
  \left\{g_1(\Gn)+
  \fr{(-1)^n i q(\Gz)}{|q(\Gn)|}
 [g_0(\Gn)+\Gr_j'] \right\}
 \fr{d\Gn}{\Gn-\Gz}.
 \label{6.16}
 \eeq
The solution  is bounded at infinity if and only if
\beq
\sum_{j=0}^{n-1}\int_{l_j}\left[\Gr_j+g(\Gx)\right]\fr{\Gx^{m-1}d\Gx}{v}=0, \quad m=1,2,\ldots, n-1.
 \label{6.17}
 \eeq
One of the constants say, $\Gr_0$,  is free. The others  have to be recovered from the 
following system  of linear algebraic equations:
  \beq
\sum_{j=1}^{n-1} (-1)^{j}I_{mj} \Gr_j =-K_m, \quad m=1,2,\ldots,n-1.
\label{6.18}
\eeq
Here,
\beq
K_m=\Gl_j\sum_{j=0}^{n-1}(-1)^j\int_{l_j^+}\fr{\Gn^{m-1}g_0(\Gn)d\Gn}{|q(\Gn)|}+I_{m0}\Gr_0.
\label{6.19}
\eeq
Since the matrix of the system (\ref{6.18}) is not singular and the coefficients $I_{mj}$ and the integrals 
$K_m$ are real, the solution to the system (\ref{6.18}) exists, unique, and real. 

When a point $\Gz$ 
traverses the loops $l_m$ the point $z=\Go(\Gz)$ traverses the contours $L_m$, 
 $$
 \Go(\Gx\pm i0)=\Gg'+\fr{c_{-1}}{\Gx-\Gz_\infty}-\fr{i\mu}{\pi\bar\tau}\left\{
 \sum_{j=0}^{n-1}(-1)^j\tilde\Gl_j\int_{l_j^+}\left[ g_1(\Gn)
 \right.\right.
$$
$$
\left.\left.
\pm (-1)^m\fr{|q(\Gx)|}{|q(\Gn)|}(g_0(\Gn)+\Gr_j')\right] \fr{d\Gn}{\Gn-\Gx}+\pi i\tilde\Gl_m[g_0(\Gx)+\Gr_m'\pm(-1)^mg_1(\Gx)]
\right\},
$$
\beq
\quad \Gz=\Gx\pm i0\in l_m^\pm, \quad m=0,1,\ldots,n-1,
\label{6.20}
\eeq 
 where $\Gg'=-i\Gg_0/\bar\tau$ is an arbitrary complex constant.
 The conformal map derived in addition to the free additive complex constant $\Gg'$ and the scaling
 parameter $c_{-1}=c_{-1}'+ic_{-1}''$ is defined up to $2n+2$ free real parameters, $a_0$, $\Gr_0$, 
 $\Gz_\infty=\Gz_\infty'+i\Gz_\infty''$, and
 $k_1,k_2,\ldots k_{2n-2}\in (-1,1)$. These parameters  have to be chosen such that
 the contours $L_m$ ($m=0,1,\ldots,n-1$) are not embedded into each other or do not intersect or touch each other.

\setcounter{equation}{0}

\section{Conclusion}

 On pursuing the goal of reconstructing 
the shapes of $n$ uniformly stressed inclusions in an unbounded elastic body subjected to antiplane uniform shear
we advanced the  technique of conformal mappings and the Riemann-Hilbert problems on a Riemann surface.
The method treats the exterior of the inclusions as the image by a conformal map of an $n$-connected slit domain with the slits lying in the real axis.
Such a map always exists for simply-, doubly-, and triply-connected domains. In the case $n\ge 4$ we considered a family
of $n$-connected domains for which their preimages lie in the same line.
The procedure requires solving consequently two inhomogeneous Schwarz problems
of the theory of analytic functions on $n$ slits or, equivalently, two Riemann-Hilbert problems
on a hyperelliptic surface of genus $n-1$. These problems were solved in terms of the Weierstrass integrals, analogs of the Cauchy integrals in a Riemann surface.
On satisfying the solvability conditions we derived a nonhomogeneous system of $n$ linear algebraic equations. The system coefficients are
basis abelian integrals, and the system matrix is nonsingular. By deriving the conformal map by quadratures and by letting a point traverse the 
double-sided slits in the parametric plane we determined the inclusions' profiles. 
The algorithm was numerically tested
for the cases $n=1$, $n=2$, and $n=3$.  
We managed to recover  families of conformal mappings  which generate sets of one, two, and three uniformly stressed inclusions.
In doubly-connected case, in addition to the complex scaling parameter $c_{-1}$ the map possesses four real parameters,
$k\in(0,1)$, $\Gz_\infty\in(-k,k)$, $a_0$, and $b_0$. In the triply-connected case it has seven real parameters,
$k\in(0,1)$, $k_1$, $k_2$   ($-k<k_1<k_2<k$), $\Gz_\infty'$, $\Gz_\infty''$, $a_0$, and $b_0$.

\vspace{2mm}

\noindent
{\bf Data accessibility.} This work does not have any experimental data.

\vspace{2mm}

\noindent {\bf Competing interests statement.}  I have no competing interests.

\vspace{2mm}

\noindent {\bf Acknowledgements.} This research originated during the author's visit to Imperial College
supported by the LMS. 

\vspace{2mm}

\noindent {\bf Funding.} This research received no specific grant from any funding agency in the public, commercial or not-for-profit sectors.

\vspace{.1in}

{\bf References}

\vspace{1mm}

\noindent 
Aksent'ev, L.A., Il'inskii, N.B.,  Nuzhin, M.T.,
Salimov, R.B. \&  Tumashev, G.G. 1980 The theory of inverse
boundary value problems for analytic functions and
its applications. (Russian) {\it Mathematical Analysis} {\bf 18}, 67-124, Akad. Nauk SSSR, Vsesoyuz. Inst. Nauchn. i Tekhn. Informatsii, Moscow.
  
\noindent 
Antipov, Y. A. \& Crowdy, D. G. 2007 Riemann-Hilbert problem for automorphic functions and the Schottky-Klein prime function. {\it Complex Anal. Oper. Theory} {\bf 1}, 317-334.

\noindent 
Antipov, Y.A. \& Silvestrov, V.V. 2007
Method of Riemann surfaces in the study of supercavitating flow around two hydrofoils in a channel,  {\it Physica D},  
{\bf 235}, 72-81.

\noindent 
Antipov, Y.A. \& Silvestrov, V.V. 2007 Method of automorphic functions in the study of flow around a stack of porous cylinders. {\it Quart. J. Mech. Appl. Math.} {\bf 60}, 337-366.

\noindent 
Antipov,  Y.A. 2017 Slit maps in the study of equal-strength cavities in 
$n$-connected elastic planar domains. {\it SIAM, J. Appl. Math.}, to appear. Preprint, 21p.:  arXiv:1703.04896v2.

\noindent 
Chaudhuri, R.A. 2003 Three-dimensional asymptotic stress field in the vicinity of the
circumferential line of intersection of an inclusion and plate
surface. {\it Int. J. Fract.} {\bf 119}, 195-222.

\noindent 
Cherepanov, G.P.  1974 Inverse problems of the plane theory of elasticity. {\it J. Appl. Math. Mech.} {\bf 38}, 
915-931. 

\noindent 
Chibrikova, L. I. \& Silvestrov, V. V. 1978 On the question of the effectiveness of the solution of Riemann's boundary value problem for automorphic functions. 
{\it Soviet Math. (Iz. VUZ)} {\bf 22}, no. 12, 85-88.

\noindent 
Courant, R.  1950 {\it Dirichlet's Principle, Conformal Mapping, and Minimal Surfaces}, New York, N.Y.:
Interscience Publishers, Inc.

\noindent  
Dai, M., Ru, C.Q. \& Gao, C.-F. 2017 Uniform strain fields inside multiple
inclusions in an elastic infinite plane under
anti-plane shear.  {\it Math. Mech. Solids}  {\bf 17} 114-128.

\noindent 
Eshelby, J.D. 1957  The determination of the elastic field of an ellipsoidal inclusion, and related problems.
{\it Proc. Roy. Soc. London A} {\bf  241}, 376-396.


\noindent 
Grabovsky, Y. \& Kohn, R.V.  1995 Microstructures minimizing the energy of a two phase elastic composite in two space dimensions. II.: The Vigdergauz microstructure. {\it J. Mech. Phys. Solids} {\bf  43},  949-972.

\noindent 
Kang, H.,  Kim, E.  \& Milton, G.W. 2008 Inclusion pairs satisfying Eshelby's uniformity property.
{\it SIAM, J. Appl. Math.} {\bf 69},  577-595.

\noindent 
Keldysh, M.V. 1939
Conformal mappings of multiply connected domains on canonical domains.
{\it Uspekhi Matem. Nauk} {\bf  6}, 90-119.

\noindent 
Liu, L.P. 2008 Solutions to the Eshelby Conjectures. {\it Proc. Roy. Soc. London A} {\bf  464}, 573-594.

\noindent 
Mityushev, V.V. \& Rogosin, S.V. 2000 {\it Constructive Methods for Linear and Nonlinear Boundary Value Problems for Analytic Functions}, London, UK: Chapman \& Hall CRC.

\noindent 
Riabuchinsky D. 1929 Sur la d\'etermination d'une surface d`apr\'es les donn\'ees
qu'elle porte. {\it C.-R. Paris} {\bf 189}  629-632.

\noindent 
Ru, C.-Q. \& Schiavone, P.  1996 On the elliptic inclusion in anti-plane shear.  {\it Math. Mech. Solids} {\bf 1},
327-333.

\noindent 
Sendeckyj,  G.P. 1970 Elastic inclusion problems in plane elastostatics. {\it Int. J. Solids Structures} {\bf 6},
1535-1543.

\noindent 
Springer, G. 1956 {\it Introduction to Riemann Surfaces}, Reading, MA:  Addison--Wesley.

\noindent 
 Vigdergauz,  S.B. 1976 Integral equation of the inverse problem of the plane theory of elasticity.
{\it J. Appl. Math. Mech.}  {\bf 40}, 518-522.

\noindent 
Wang, X. 2012 Uniform fields inside two
non-elliptical inclusions. {\it Math. Mech. Solids}  {\bf 17} 736-761.

\end{document}